\documentclass[11pt]{article} 

\usepackage{times}
\usepackage[square]{natbib}
\usepackage{epsfig,amsmath,amssymb,cite}
\usepackage{graphicx}

\newtheorem{teor}{Theorem}[section]
\newtheorem{lema}{Lemma}[section]

\newtheorem{probld}{Open Problem}
\newtheorem{cnj}{Conjecture}
\newtheorem{exempl}{Example}[section]

\def\bull{\vrule height 1.4ex width .8ex depth -.1ex }
\newenvironment{dem}{\noindent{\bf Proof}\ \ }{\bull}

\def\sqr#1#2{{\vcenter{\hrule height.#2pt
    \hbox{\vrule width.#2pt height#1pt \kern#1pt
    \vrule width.#2pt}
    \hrule height.#2pt}}}

\def\a{\alpha}
\def\b{\beta}

\def\settabs{\setbox\tabs=\null \futurelet\next\sett@b}

\def\newblock{\hskip .11em plus.33em minus.07em}

\title{About the uncountability of the number of irrational\\ 
powers of irrational numbers evaluated as rationals\\
and solutions' estimation for $x^x=y$ and $x^{x^x}=y$}
\author{Anca Andrei\\
Department of Mathematics\\
University of Texas at Austin\\
Austin, TX, U.S.A.\\
{\normalsize \tt{anca.andrei@utexas.edu}}\\
}

\begin{document}
\maketitle

\begin{abstract}
It is well known that numbers are often used to define numbers. 
For example, two integers are used to define a rational number and two real numbers are used to
define a complex number. 
It might be expected that an irrational power of an irrational
number would be an irrational number. 
Despite this expectation, it is possible for an irrational
power of an irrational number to be a rational number. For many generations, this circulated
as a non-constructive proof by contradiction in logic for discrete mathematics
textbooks and college courses \citep{Ros2011}. 
Since the '80s, a constructive proof circulated orally, such as 
$(\sqrt{2})^{log_{\sqrt{2}} 3}$ equals to $3.$ A written proof was published
in 2008 by Lord \citep{Lor2008}.

The first contribution of this paper is to show that there is an
uncountable number of such pairs of irrational numbers such that
the power of one to the other is a rational number.

Marshall and Tan answered the question of
whether there is a single irrational number $a$ such that $a^a$ is 
rational \citep{MaT2012}. They proved that given $I=((\frac{1}{e})^{\frac{1}{e}}, \infty)$,
then every rational number in $I$ is either of the form $a^a$ for an 
irrational $a$ or is in the very thin set $\{$1, 4, 27, 256, ..., $n^n$, ...
$\}.$ 
It seems a challenging task to analytically solve the equation
$x^x=y$ for any real $y$. To the best of our knowledge, there is no work
on finding $x$ from a given $y.$ We proved that $ln(ln(y))<x$ for $y>e$ and
$x<ln(y),$ for $y>e^e.$ 

Hence, the second contribution of this paper is to estimate the real $x$ such that
either one of the equations $x^x=y$ or $x^{x^x}=y$ holds, for a given $y.$
\end{abstract}









\section{Introduction}

This work establishes that there exists an 
uncountable number of irrational numbers $a$ and $b$ such that
$a^b$ is a rational number.

In 1900, David Hilbert posed a list of challenging problems,
including a general problem referring to whether $2^{\sqrt{2}}$ 
and ${\sqrt{2}}^{\sqrt{2}}$ are transcendental numbers, among other
interesting powers. 
A number is {\em algebraic} 
if and only if it is the root of a polynomial equation with rational 
coefficients. On the contrary, a real number is {\em transcendental} 
if and only if the number is not algebraic. In 1934, 
Gelfond and Schneider proved, independently 
\citep{Gel34a,Gel34b,Sch34a,Sch34b}, that if 
$\alpha$ and $\beta$ are algebraic numbers with 
$\alpha \neq 0,$ $\alpha \neq 1,$ and $\beta \notin \mathbb{Q}$,
then $\alpha^\beta$ is transcendental. Following Gelfond-Schneider 
Theorem, both numbers $2^{\sqrt{2}}$ and ${\sqrt{2}}^{\sqrt{2}}$ 
are transcendental.

In many discrete mathematics books, such as \citep{Ros2011}, 
a non-constructive proof 
for the existence of irrational powers of irrational numbers that 
are rational is provided. For example, the proof described in 
\citep{Ros2011} does not indicate any such numbers, but only a 
non-constructive proof based on the fact that 
${\left((\sqrt{2})^{\sqrt{2}}\right)}^{\sqrt{2}}=2$.

We continue these works by stating a few results about the 
irrational powers of irrational numbers which are evaluated as
rationals.
It is easy to prove that ${\left((\sqrt{n})^{\sqrt{n}}\right)}^{\sqrt{n}}$
is rational if and only if $n$ is an even positive number, since 
${\left((\sqrt{n})^{\sqrt{n}}\right)}^{\sqrt{n}}=
n^{\frac{n}{2}}$. For example, 
${\left((\sqrt{4})^{\sqrt{4}}\right)}^{\sqrt{4}}=
4^{\frac{4}{2}}=16$, hence rational, while 
${\left((\sqrt{5})^{\sqrt{5}}\right)}^{\sqrt{5}}=
25{\sqrt{5}}$, which is irrational.

In 2008, Nick Lord \citep{Lor2008} published a constructive proof 
based on the fact that $(\sqrt{2})^{log_{\sqrt{2}} 3}$ $= 3$. 
Both $\sqrt{2}$ and ${log_{\sqrt{2}} 3}$ are irrational numbers, 
but their power composition $(\sqrt{2})^{log_{\sqrt{2}} 3}$ is a 
rational number. Lord generalizes this problem to more than just one 
number. He proved that if $\frac{m}{n}$ is any positive rational
with $m\neq n$, and $p$ is any prime 
which is neither a factor of $m$ nor of $n$, then 
${\sqrt{p}}$ and ${log_{\sqrt{p}} \frac{m}{n}}$ are both irrational.  
Furthermore, their power composition $(\sqrt{p})^{log_{\sqrt{p}} \frac{m}{n}}$ 
is the rational $\frac{m}{n}$.

An interesting question is how many irrationals ${\a}$
and ${\b}$ exist such that ${\a^\b}$ is rational. Lord's result implies
that there are (at least) a countable set of numbers ${\a}$ and ${\b}$ 
such that ${\a^\b}$ is a rational. A {\em countable} set is a set with the the same
cardinality as the set of natural numbers. A set that is
not countable is called {\em uncountable}. Cantor proved that the rational set
of numbers and the algebraic set of numbers are both countable sets. He also proved
that the set of real numbers is uncountable. 

Since the set of algebraic numbers is countable and the set of real
numbers is uncountable, it follows that the set of transcendental numbers 
is uncountable. 

Marshall and Tan answered the question of
whether there is a single irrational number $a$ such that $a^a$ is 
rational \citep{MaT2012}. They proved that given $I=((\frac{1}{e})^{\frac{1}{e}}, \infty)$,
then every rational number in $I$ is either of the form $a^a$ for an 
irrational $a$ or is in the very thin set $\{$1, 4, 27, 256, ..., $n^n$, ...
$\}.$ 
It seems a daunting task to analytically solve the equation
$x^x=y$ for any real $y$. To the best of our knowledge, there is no work
on finding $x$ from a given $y.$ We proved that $ln(ln(y))<x$ for $y>e$ and
$x<ln(y),$ for $y>e^e.$ 

The paper's contribution is two-fold:
\begin{enumerate}
\item  to show that there is an
uncountable number of pairs of irrational numbers such that
the power of one to the other is a rational number;
\item to estimate the real $x$ such that
either one of the equations $x^x=y$ or $x^{x^x}=y$ holds, for a given $y.$
\end{enumerate}


\section{The result about uncountability of irrationals power of irrational
that are rational}

\noindent
\begin{lema}\label{l.1.}
If $a$ is a positive transcendental number and $c$ is a positive integer, 
then $log_a c$ is an irrational number.
\end{lema}

\begin{dem}
By contradiction, we assume that $log_a c$ is a rational number. 
Hence, there exists a positive integer $p$ and a non-zero integer $q$ 
such that $log_a c = \frac{p}{q}$. This means $a^\frac{p}{q} = c$ which 
implies $a = (c)^\frac{q}{p}$. The number $(c)^\frac{q}{p}$ is
algebraic because it is the root of the polynomial $X^p - c^q =0$
with rational 
coefficients. On the other hand, $a$ is a transcendental number, so the
equality $a = (c)^\frac{q}{p}$ represents a contradiction. 
As a consequence, $log_a c$ is an irrational number and
$a^{log_a c}$ is rational. \hfill\bull
\end{dem}
\vskip 5pt

\noindent
\begin{lema}\label{l.2.}
If $a$ is a positive transcendental number, then there exists an 
irrational number $b$ such that $a^b$ is a rational number.
\end{lema}

\begin{dem}
Let $a$ be a positive transcendental number.
Let $b= log_a c$, where $c$ is a positive rational 
integer. Obviously, $log_a c$ is well defined because $a$ and $c$ are 
positive numbers. Since $a$ is transcendental, $a\neq 1$
so the logarithmic function $log_a$ is well defined.

By Lemma \ref{l.1.}, $log_a c$ is an irrational number.
Hence $a^b = a^{log_a c}$ is rational. 
\hfill \bull
\end{dem}
\vskip 5pt

We are ready to prove the main result of this section.
\vskip 5pt

\noindent
\begin{teor}\label{t.1.}
There exists an uncountable set of irrational numbers $a$ 
and $b$ such that $a^b$ is a rational number. 
\end{teor}

\begin{dem}
We choose ${a}$ to be a positive transcendental number and 
$b= log_a c$, where $c$ is an arbitrary positive rational number.
According to Lemma \ref{l.1.}, it follows that ${b}$ is an irrational number. 
Based on Cantor's Theorem, the set of positive transcendental numbers is 
uncountable. 
Combining these aforementioned statements, this theorem 
is proved.
\hfill \bull
\end{dem}
\vskip 5pt

\section{An estimation of the solutions of $x^x=y$ and $x^{x^x}=y$, for a given $y$}\label{s.3.}

To obtain an estimation of the solution of $x^x=y$, we need an 
intermediate support lemma.

\begin{lema}\label{l.3.}
If $z$ is a positive real, then $e^z>z+ln(z).$
\end{lema}

\begin{dem}
If $z\in(0,1]$, then $e^z>1>z+ln(z)$ and the conclusion holds.

If $z\in(1,\infty)$, then we define $f$: $\mathbb{R}\to\mathbb{R}$ given by $f(z)=e^z-z-ln(z).$ 
Obviously, $f$ is a continuous and differentiable function on $(1, \infty)$, with
$\displaystyle{f'(z)=\frac{z\cdot e^z-z-1}{z}}.$
Defining $g$: $\mathbb{R}\to\mathbb{R}$ by $g(z)=z\cdot e^z-z-1,$ we see that function 
$g$ is continuous and has its first continuous derivatives.
Its first derivative is $g'(z)=e^z+z\cdot e^z-1$ and its second derivative is 
$g''(z)=(z+2)\cdot e^z.$ Clearly, $g''(z)>0$, for $z>1.$ Hence $g'$ is a 
monotone increasing function, so $g'(z)>g'(1).$ Since $g'(1)=2\cdot e-1$ is a positive real,
$g$ is also a monotone increasing function. Hence $g(z)>g(1).$
Since $g(1)=e-2$ is a positive real, $g(z)>0$. Thus $f'(z)>0,$ $\forall z>1.$
This implies that $f$ is a monotone increasing function.
Hence $f(z)>f(1)=e-1$. Since this quantity is positive, the inequality
of this lemma holds for $z>1$, too.\hfill\bull
\end{dem}

Considering the equation $x^x=y$ for any real given $y$,
Theorem \ref{t.2.} estimates the solution $x$.

\begin{teor}\label{t.2.}
Let us consider the equation $x^x=y$, for a given real
$y$. Then, the following estimations of solution $x$ hold:
\begin{enumerate}
\item[(a)] if $y>e^e$, then $x<ln(y)$.
\item[(b)] if $y>e$, then $x>ln(ln(y))$.
\end{enumerate}
\end{teor}

\begin{dem}
(a) Since $y=x^x$ and $y>e^e$, the inequality $x<ln(y)$ is equivalent to $x<ln(x^x)$,
which is equivalent in turn to $x<x\cdot ln(x)$. Since $x^x=y>e^e$, it follows that $x>e$. 
Hence it follows that $ln(x)>1$, so $x<x\cdot ln(x)$ holds.

(b) By substituting $y=x^{x^x}$, we get an equivalent inequality
$x>ln(ln(x^x))$, where $x\cdot ln(x)>1$. This inequality is equivalent to 
$x>ln(x\cdot ln(x))$, which is in turn equivalent to
$x>ln(x)+ln(ln(x))$ for any real $x$ such that 
$x\cdot ln(x)>1$. Since $x\cdot ln(x)>1$, it means that $x>1$. Thus, there
exists $z>0$ such that $x=e^z$. The inequality to be proved becomes:
$e^z>z+ln(z)$, 
Since this is the inequality from Lemma \ref{l.3.}, 
we conclude the proof.
\hfill\bull
\end{dem}

We are now ready to handle the exponentiated version of the previous
equation, that is, $x^{x^x}=y$, for a given $y$.

\begin{lema}\label{l.4.}
If $x>e$, then the following statements hold:
\begin{enumerate}
\item[(a)] $x<x\cdot ln(x)+ln(ln(x));$
\item[(b)] $x>ln(x\cdot ln(x)+ln(ln(x))).$
\end{enumerate}
\end{lema}

\begin{dem}
(a) Since $x>e$, it means $ln(x)>1,$ so $ln(ln(x))>0.$ Hence $x<x\cdot ln(x)+ln(ln(x));$

(b) We prove a stronger inequality, that is, $x>ln(x+1)+ln(ln(x)).$
Since $x>e$, it means $\exists z>1$ such that $x=e^z.$ We need to
prove $e^z>ln(z\cdot(e^z+1)).$ 
We define $f$: $\mathbb{R}\to\mathbb{R}$ given by $f(z)=e^z-ln(z\cdot(e^z+1)).$
Obviously, we observe that $f$ is continuous and differentiable with first derivative
$f'(z)=\frac{z\cdot(e^z)^2-e^z-1}{z\cdot (e^z+1)}.$
So, $f(z)>f(1)>0.$ Hence $x>ln(x+1)+ln(ln(x)).$
The right hand side can be rewritten as:
$$ln(x+1)+ln(ln(x))=ln((x+1)\cdot ln(x))= ln(x\cdot ln(x)+ ln(x)).$$ 
But $ln(x)<x,$ so $ln(x)>ln(ln(x))$ by substituting $x$ with $ln(x).$ 
Hence $$ln(x\cdot ln(x)+ ln(x))>ln(x\cdot ln(x)+ln(ln(x))).$$ 
Therefore $x>ln(x\cdot ln(x)+ln(ln(x))).$ The lemma is proved. \hfill\bull
\end{dem}
\vskip 5pt

\begin{teor}\label{t.3.}
Let us consider the equation $x^{x^x}=y$, for a given $y>e^{e^e}.$ Then
the following estimations of solution $x$ hold:
\begin{enumerate}
\item[(a)] $x<ln(ln(y));$
\item[(b)] $x>ln(ln(ln(y))).$
\end{enumerate}
\end{teor}

\begin{dem}
Since $x^{x^x}=y$ and $y>e^{e^e},$ it means that $x>e$.

(a) The inequality $x<ln(ln(y))$ is equivalent to 
$x<ln(x^x\cdot ln(x))$ or $x<x\cdot ln(x)+ln(ln(x)).$
By item (a) of Lemma \ref{l.4.}, this inequality holds.

(b) Similarly, $x>ln(ln(ln(y)))$ is equivalent to 
$x>ln(x\cdot ln(x)+ln(ln(x))).$ According to item (b) of Lemma \ref{l.4.},
this inequality holds. \hfill\bull
\end{dem}
\vskip 5pt


\noindent
{\bf Acknowledgments.} We thank Dr. Ted W. Mahavier from Lamar 
University and to Dr. Ronny Hadani from University of Texas at Austin 
for their valuable comments while proofreading this work.

\section*{Conclusions}\label{s.4.}

This paper proved that there exists an uncountable number of 
irrational numbers $a$ and $b$ such that $a^b$ is a rational 
number. We also estimate the solutions x for both equations
$x^x=y$ and $x^{x^x}=y.$ 
  Our results represent extensions of the works done by
Lord in \citep{Lor2008} and Marshall and Tan in \citep{MaT2012}.


\end{document}